\documentclass[12pt,a4paper]{article}
\usepackage{amssymb,amsmath}
\usepackage{amsfonts}
\usepackage{centernot}
\usepackage{mathtools}
\usepackage{stmaryrd}
\usepackage{a4}
\newcommand{\para}{\par\vspace{.25cm}}

\newtheorem{prop}{Proposition}
\newtheorem{theorem}{Theorem}

\newtheorem{cor}{Corollary}
\newtheorem{remark}{Remark}
\usepackage{multirow}

\begin{document}
\baselineskip 18pt

\title{\bf Integral group rings with all central units trivial }
\author{Gurmeet K. Bakshi  \\ {\em \small Centre for Advanced Study in
Mathematics,}\\
{\em \small Panjab University, Chandigarh 160014, India.}\\{\em
\small email: gkbakshi@pu.ac.in} \and  Sugandha Maheshwary {\footnote {Research supported by IISER, Mohali, India, is gratefully acknowledged} \footnote{Corresponding author}}
\\ {\em \small Indian Institute of Science Education and Research, Mohali,}\\
{\em \small Sector 81, Mohali (Punjab)-140306, India.}
\\{\em \small email: sugandha@iisermohali.ac.in}
\and Inder Bir S. Passi \\  {\em \small Centre for Advanced Study
in Mathematics,}\\ {\em \small Panjab University, Chandigarh-160014, India} \\
{\small \& }\\
{\em \small Indian Institute of Science Education and Research, Mohali,}\\
{\em \small Sector 81, Mohali (Punjab)-140306, India.}\\{\em \small email: ibspassi@yahoo.co.in } }
\date{}
{\maketitle}

\begin{abstract}\noindent {}
The object of study in this paper is the finite groups whose integral group rings have only trivial central units. Prime-power groups and metacyclic groups with this property are characterized. Metacyclic groups are classified according to the central height of the unit groups of their integral group rings.
\end{abstract}\vspace{.25cm}
 \vspace{.25cm} 

\section{Introduction}Given a finite group $G$, let $\mathcal{U}(\mathbb{Z}[G])$ be the group of units of the integral group ring $\mathbb{Z}[G]$ and let $\mathcal{Z}(\mathcal{U}(\mathbb{Z}[G]))$ be its centre. The structure of $\mathcal{Z}(\mathcal{U}(\mathbb{Z}[G]))$ has been a \linebreak subject of intensive research (see e.g. \cite{BM2,Ferraz,good,del,jesper,jesp3,jesp2,parm,Li,sehg}). Clearly $\mathcal{Z}(\mathcal{U}(\mathbb{Z}[G]))$ contains the so-called \emph{trivial central units} $\pm g$, $g \in \mathcal{Z}(G)$, the centre of~$G$. Thus there arises the problem of characterizing the groups $G$ having the property that all central units of the integral group ring $\mathbb{Z}[G]$ are trivial. For notational convenience, we refer to such groups as groups with the {\sf cut}-property. \emph{Abelian groups with the  {\sf cut}-property are precisely those having exponent $1,\, 2, \,3,\, 4$ or ~$6$} \cite{Hig}. It has been shown by Ritter and Sehgal \cite{sehg} that a finite group $G$ has the {\sf cut}-property if, and only if, for all $x\in G$ and for every natural number $j$, relatively prime to the order of~$G$,
\begin{equation}\label{Eq1}
    x^{j}~is~ conjugate~to ~either ~ x ~or~x^{-1}.
\end{equation}The above characterization is equivalent (see \cite{sehg} \& \cite{hupp}, p.\,545) to saying that, if $\mathbb{Q}[G]\cong \bigoplus_{i}M_{n_{i}}(\mathbb{D}_{i})$, is the Wedderburn decomposition of the rational group algebra $\mathbb{Q}[G]$, where $M_{n_{i}}(\mathbb{D}_{i})$ denotes the algebra of $n_{i}\times n_{i}$ matrices over division ring $\mathbb{D}_{i}$, then
  \begin{equation}\label{Eq2}
   ~the ~centre~ \mathcal{Z}(\mathbb{D}_{i})~ of~ each ~\mathbb{D}_{i}~ is~ either~ \mathbb{Q} ~or~ ~\mathbb{Q}(\sqrt{-d}),~d\in \mathbb{N}.
\end{equation}
In particular, \emph{the  {\sf cut}-property is quotient closed}, i.e., if $G$ is a group with the \linebreak{\sf cut}-property, then so is every quotient of $G$. Our aim in this paper is to \linebreak investigate further the {\sf cut}-property for finite groups. The study of this property \linebreak reveals, as expected, that it is a fairly restrictive property. After making some \linebreak general observations on the groups with the {\sf cut}-property, we examine this\linebreak property in detail  for nilpotent groups and metacyclic groups.\\

In Section 2, we prove (Theorem \ref{Th1}) that the order of a
group with the\linebreak {\sf cut}-property must necessarily be divisible by 2 or 3. Consequently, since this property is quotient closed, a nilpotent group with the  {\sf cut}-property must
be a $(2, \,3)$-group. We note that the above-mentioned  characterization (\ref{Eq1}) of  groups with the {\sf cut}-property  can be considerably simplified in the case of $2$-groups and $3$-groups. We show that a $3$-group $G$ has the  {\sf cut}-property if, and only if, for all $x \in G$, $x^{2}$ is conjugate to $x^{-1}$ (Theorem~\ref{Th3}) and a $2$-group $G$ has the {\sf cut}-property if, and only if, for each $x \in G$, $x^{3}$ is conjugate to either $x$ or $x^{-1}$ (Theorem \ref{Th2}). It is interesting to note that, for $3$-groups, the {\sf cut}-property is closed under direct sum, but it is not so  for $2$-groups.\\

In Section 3, as a contribution to understand the finite solvable groups with the {\sf cut}-property, we examine metacyclic groups, and give a complete classification of such groups with this property. It turns out that, up to isomorphism, only finitely many metacyclic groups have the {\sf cut}-property, of which we give a complete list (Theorem \ref{Th7}).\\

It is known \cite{satya} that the central height of $V(\mathbb{Z}[G])$, the group of units of augmentation $1$ in $\mathbb{Z}[G]$, is at most $2$. We compute precisely the central height of $V(\mathbb{Z}[G])$, when $G$ is a metacyclic group (Theorem \ref{Th6}).

\section{Finite groups with $\mathcal{Z}(\mathcal{U}(\mathbb{Z}[G]))=\pm \mathcal{Z}(G)$}

We begin by observing that the {\sf cut}-property has strong bearing on the centre of the group. If $G$ has the {\sf cut}-property, then so does its centre  $\mathcal{Z}(G)$. For if $\mathcal{Z}=\mathcal{Z}(G)$, then $\mathcal{U}(\mathbb{Z}[\mathcal{Z}])\subseteq \mathcal{Z}(\mathcal{U}(\mathbb{Z}[G]))$. In case  $G$ has the {\sf cut}-property, $\mathcal{Z}(\mathcal{U}(\mathbb{Z}[G])) =\pm\mathcal{Z}(G) $, and consequently,  in that case, $\mathcal{U}(\mathbb{Z}[\mathcal{Z}]) =\pm\mathcal{Z}$.\\

\noindent For a finite group $G$, let   $\pi(G)$ denote  the set of primes which divide the order of~$G$.
 \begin{theorem}\label{Th1}
If $G$ has the {\sf cut}-property, then either $2$ or $3 \in \pi(G)$. Furthermore, if $G$ is nilpotent, then $\pi(G)\subseteq \{2,\,3\}$.
\end{theorem}
{\bf Proof.} Let $G$ be a group with the {\sf cut}-property.  If $2 \not \in \pi(G)$, then, by the \linebreak odd-order theorem \cite{Feit},  $G$ is solvable  and hence $G\neq G'$, the derived subgroup of $G$. Let $p$ be a prime divisor of $[G:G']$, the index of $G'$ in $G$. Since the  \linebreak  {\sf cut}-property is quotient-closed, $G/G'$ has the {\sf cut}-property. Thus it follows that $G/G'$ is a $3$-group, and hence $3\in\pi(G)$. Further, if $G$ has the {\sf cut}-property, then so does  $G/\mathcal{Z}(G)$. Thus, if $G$ is nilpotent group, then each quotient in the upper central series of $G$, being abelian with the {\sf cut}-property, is of exponent $2,3,4$ or $6$ and hence $\pi(G)\subseteq \{2,\,3\}$.  $\Box$\\

We next investigate the $2$-groups and $3$-groups with the {\sf cut}-property.

\begin{theorem}\label{Th3}
A finite $3$-group has the {\sf cut}-property if, and only if,  $x^{2}$ is conjugate to $ x^{-1}$, for all $ x~\in G$. Furthermore, if a $3$-group $G$ has the {\sf cut}-property, then each quotient in the lower, as well as the upper, central series of $G$ is of exponent $3$.
 \end{theorem}
{\bf Proof.} In view of the characterization (\ref{Eq1}) and the fact that $2$ is a primitive root modulo $3^{n}$, for all $n\geq 1$, a $3$-group $G$ has the {\sf cut}-property if, and only if, $x^{2}$ is conjugate to $ x$ or $x^{-1}$, for all $ x~\in G$.
\par\vspace{.25cm}
Let $G$ be a $3$-group and let $$G=\gamma_{1}(G)\supset \gamma_{2}(G) \supset ...\supset \gamma_{n}(G)=\langle 1 \rangle$$ be its lower central series.
We claim that $x^2$ is not conjugate to $x$ for any $1\neq x\in G$.  For, given $1\neq x\in G$, there exists $i$ $(1\leq i<n)$ such  that $x\in \gamma_{i}(G)\setminus\gamma_{i+1}(G)$.  If $x^{2}$ were conjugate to $x$, then $x^{2}=g^{-1}xg$ for some $g\in G$, and therefore $x=x^{-1}g^{-1}xg \in \gamma_{i+1}(G)$, which contradicts the assumption. This proves the first statement. \par\vspace{.25cm}
Thus, if $G$ is a $3$-group with the {\sf cut}-property, then for every $1\neq x\in G$ and $1\leq j < n$, there exits $y\in G$ such that $x^{2}=y^{-1}x^{-1}y$ and therefore, $x^{3}\in \gamma_{j+1}(G)$.
\par\vspace{.25cm}
Next, let $$\langle1 \rangle=\mathcal{Z}_{0}(G)\subset \mathcal{Z}_{1}(G) \subset ...\subset \mathcal{Z}_{n}(G)=G$$ be the upper central series  of $G$. Then, $\mathcal{Z}_{1}(G)=\mathcal{Z}(G)$ is of exponent $3$ and as the {\sf cut}-property is quotient closed, it follows inductively that, for every $i$, $\mathcal{Z}_{i+1}(G)/\mathcal{Z}_{i}(G)$ is of exponent $3$. $\Box$
\begin{cor}\label{Co3}The {\sf cut}-property is direct sum closed for $3$-groups.
\end{cor}
{\bf Proof.} Let $H$ and $K$ be $3$-groups with the {\sf cut}-property and let $G= H \times K$ be the direct sum of $H$ and $K$. Further, let $g=(h,\,k)\in G$, $h \in H$, $k\in K $. Then, $g^{2}=(h^{2},\,k^{2})$ is conjugate to $(h^{-1},\ k^{-1})=g^{-1}$ and hence, by Theorem \ref{Th3}, $G$ has the {\sf cut}-property. $\Box$
\par\vspace{.2cm}
In view of Eq.\,(\ref{Eq1}) and the fact that $\mathcal{U}(\mathbb{Z}/2^{n}\mathbb{Z})=\langle \pm 1\rangle \oplus \langle 3 \rangle$, we have the following analogous result for $2$-groups.
\begin{theorem}\label{Th2}
A finite $2$-group has the {\sf cut}-property if, and only if,  $x^{3}$ is conjugate to $ x$ or $x^{-1}$, for all $ x~\in G$. Furthermore, if a $2$-group $G$ has the {\sf cut}-property, then each quotient in its lower central,  as well as the upper central, series is of exponent $2$ or $4$.

\end{theorem}

\begin{remark}
Direct sum of $2$-groups with the {\sf cut}-property does not necessarily  have the {\sf cut}-property. For example, let\par \centerline{$H=\langle a,~b~|~ a^{8}=b^{2}=1,~b^{-1}ab=a^{3}\rangle$ and  $K=\langle x~|~ x^{4}=1\rangle.$}\par\noindent It is easy to check that both $H$ and $K$ have the {\sf cut}-property. However, $g=(a,\,x)\in H \times K$, is such that $g^{3}=(a^{3},\,x^{3})$ is neither  conjugate to $g$ nor to  $g^{-1}$.
\end{remark}

 An interesting class of $p$-groups with the {\sf cut}-property is provided by Camina groups. Recall that a group $G$ is called a \emph{Camina group}, if $G' \neq G$ and, for every $g \not\in G'$, the coset $gG'$ is a conjugacy class. For a normal subgroup $N$ of $G$, the pair $(G,N)$ is said to be a \emph{Camina pair}, if each element $g \not\in N$ is conjugate to all the elements in the coset $gN$, so that $G$ is a Camina group, if $(G,G')$ is a Camina pair. Camina $p$-groups have nilpotency class at most three, and the factors of their lower central series are of exponent $p$ \cite{dark,MacD}.

\begin{theorem}\label{p13} Non-abelian Camina $p$-groups have the {\sf cut}-property for  $p=2,\ 3$.
\end{theorem} {\bf Proof.}  Let $G$ be a Camina $3$-group of class 2. Then $G/G'$ and $G'$ are both of exponent 3. It thus immediately follows that  $x^2$ is conjugate to $x^{-1}$ for every   $x\in G$, and consequently, by Theorem \ref{Th3}, $G$ has the {\sf cut}-property.
\par
In case $G$ is a Camina $2$-group of class 2, then $G$ has exponent 4, and hence $x^3=x^{-1}$ for every $x\in G$. Thus, in view of Theorem \ref{Th2}, $G$ has the {\sf cut}-property.  \par
Next, let $G$ be a $p$-group of nilpotency class $3$, $p=2,\ 3$.  Then $G/\gamma_{3}(G)$ is a Camina $p$-group of class $2$, and so it has the {\sf cut}-property. Furthermore, $(G,\gamma_{3}(G))$ forms a Camina pair (\cite{MacD}, Theorem 5.2), which ensures, by Theorems \ref{Th3} \& \ref{Th2}, that $G$ has the {\sf cut}-property. $\Box$

\section{Finite metacyclic groups}
\subsection{Metacyclic groups  with the {\sf cut}-property}
In this section, we determine the metacyclic groups with the {\sf cut}-property.\\
We begin by recalling some results from \cite{Oli} .\\

\vspace{-0.25cm} Let $G$ be a finite group and $K$ be a normal subgroup of subgroup $H$ of $G$. Define $\hat{H}:=\frac{1}{|H|} \sum_{h \in H} h$ and $$\varepsilon (H, K) :=  \begin{cases} \hat{H}, & H = K; \\ \prod (\hat{K}-\hat{L}), &  {\rm otherwise,} \end{cases} $$ where $|H|$ denotes the order of $H$ and $L$ runs over the normal subgroups of $H$, which are minimal over the normal subgroups of $H$ containing $K$ properly.\\

\noindent A {\it strong Shoda pair} of $G$ is a pair $(H,K)$ of subgroups of $G$ with the property that \begin{description}
\item[(i)] $K$ is normal in $H$ and $H$ is normal in $N_{G}(K)$, the normalizer of $K$ in $G$;
 \item[(ii)] $H/K$ is cyclic and a maximal abelian subgroup of $N_{G}(K)/K$;
 \item[(iii)] the distinct $G$-conjugates of $\varepsilon(H,K)$ \index{$\varepsilon(H,K)$} are mutually orthogonal.
     \end{description}
It is known that if $(H,K)$ is a strong Shoda pair of $G$, then $e(G,H,K)$, the sum of distinct $G$-conjugates of $\varepsilon(H,K)$, is a primitive central idempotent of $\mathbb{Q}[G]$ (\cite{Oli}, Proposition 3.3).\\

\noindent Furthermore, a strong Shoda pair $(H, K)$ of $G$ with $H$ normal in $G$ has an easier description in view of the following:

\begin{prop}\label{p13}{\rm(\cite{Oli}, Proposition 3.6)} Let $(H,K)$ be a pair of subgroups of a finite group $G$ satisfying the following conditions:
\begin{description}
  \item[(i)] $K$ is normal in $H$ and $H$ is normal in $G$,
  \item[(ii)] $H/K$ is cyclic and a maximal abelian subgroup of $N_{G}(K)/K$.
\end{description}

\noindent Then $(H,K)$ is a strong Shoda pair and hence $e(G,H,K)$ is a primitive central idempotent
of $\mathbb{Q}[G]$.

 \end{prop}
\para  A big advantage of determining the primitive central idempotent of $\mathbb{Q}[G]$, \linebreak using strong Shoda pair $(H,K)$ of $G$ is that one can describe the structure of the simple component $\mathbb{Q}[G]e(G,H,K)$ of $\mathbb{Q}[G]$, corresponding to the primitive central idempotent $e(G,H,K)$, in view of the next proposition.\\

Note that $R*^{\sigma}_{\tau}G$ denotes the crossed product of group $G$ over the coefficient ring $R$ with action $\sigma$ and twisting $\tau$ \cite{Pass}. Recall that a classical crossed product is a crossed product $K*^{\sigma}_{\tau}G$, where $K/F$ is a finite Galois extension, $G=\operatorname{Gal}(K/F)$ is the Galois group of $K/F$ and $\sigma $ is the natural action of $G$ on $K$. The classical crossed product is denoted by $(K/F,\tau)$ \cite{rei}.\\
\vspace{-0.25cm}
\begin{prop}\label{p11}{\rm(\cite{Oli}, Proposition 3.4)} Let $(H,K)$ be a strong Shoda pair of a finite group $G$ with $k=[H:K]$, $N=N_{G}(K)$, $n=[G:N]$, $x$ a generator of $H/K$ and $\phi: N/H \mapsto N/K$ a left inverse of the projection $N/K \mapsto N/H$. Then,
$$ \mathbb{Q}[G]e(G,H,K)\cong M_{n}(\mathbb{Q}(\zeta_{k})*^{\sigma}_{\tau}N/H),$$
where $\zeta_{k}$ denotes the $k^{th}$ root of unity, the action $\sigma$ and the twisting $\tau$ are given by $\zeta_{k}^{\sigma(a)}=\zeta_{k}^{i}$, if $x^{\phi(a)}=x^{i};$ $\tau(a,b)=\zeta_{k}^{j}$, if $\phi(ab)^{-1}\phi(a)\phi(b)=x^{j}$, for $a,b \in N/H$ and integers $i$ and $j$.
   \end{prop}

 \begin{remark} The action $\sigma$ in the above proposition is faithful and therefore \linebreak $\mathbb{Q}(\zeta_{k})*^{\sigma}_{\tau}N/H$ can be described as a classical crossed product $(\mathbb{Q}(\zeta_{k})/F,\tau)$, where $F$ is the centre of the algebra and $\sigma$ is determined by the action of $N/H$ on $H/K$. In this way, the Galois group $\operatorname{Gal}(\mathbb{Q}(\zeta_{k})/F)$ can be identified with $N/H$ and with this identification, $F=\mathbb{Q}(\zeta_{k})^{N/H}$, the fixed field of $N/H$.   \end{remark}

 We now state the main theorem of this section, which classifies all metacyclic groups with the {\sf cut}-property.\\
\vspace{-0.25cm}

\begin{theorem}\label{Th7} Let $G$ be a finite non-abelian metacyclic group given by \begin{equation}\label{Eq3} G=\langle a,~b~|~ a^{n}=1,~b^{t}=a^{\ell},~b^{-1}ab = a^ {r}\rangle ,\end{equation} where $n,\,t,\,r,\,\ell$ are natural numbers such that $r^{t}\equiv 1~({\rm mod}~n)$, $\ell r \equiv \ell~({\rm mod}~n)$ and $\ell$ divides $n$. Then, $G$ has the {\sf cut}-property if, and only if, $G$ is isomorphic to one of the following:
 \begin{quote} $ \langle a,~b~|~ a^{n}=1,~b^{t}=1,~b^{-1}ab = a^ {n-1}\rangle,~~t=2,4,6,~~n=3,4,6;$ \\
 $\langle a,~b~|~ a^{4}=1,~b^{t}=a^{2},~b^{-1}ab = a^ {3}\rangle,~t=2,4,6$; \\
  $ \langle a,~b~|~ a^{6}=1,~b^{2}=a^{3},~b^{-1}ab = a^ {5}\rangle;$ \\
  $\langle a,~b~|~ a^{n}=1,~b^{\varphi(n)}=1,~b^{-1}ab = a^ {\lambda_{n}}\rangle,~n=5,7,9,10,14,18$;\\
  $\langle a,~b~| ~ a^{n}=1,~b^{\frac{\varphi(n)}{j}}=1,~b^{-1}ab = a^ {\lambda_{n}^{2}}\rangle,~~j=1,2,~~n=7,9,14,18$;\\
  $\langle a,~b~|~ a^{8}=1,~b^{t}=1,~b^{-1}ab = a^ {r}\rangle,~t=2,4,~r=3,5$;\\
   $\langle a,~b~|~ a^{12}=1,~b^{t}=1,~b^{-1}ab = a^ {5}\rangle,~t=2,4$;\\
  $\langle a,~b~|~ a^{12}=1,~b^{t}=a^{\ell},~b^{-1}ab = a^ {7}\rangle,~t=2,6,~\ell=t,12$;\\
  $\langle a,~b~| ~a^{15}=1,~b^{4}=1,~b^{-1}ab = a^ {2}\rangle$;\\
  $\langle a,~b~|~ a^{16}=1,~b^{4}=1,~b^{-1}ab = a^ {r}\rangle$,~$r=3,5$;\\
  $\langle a,~b~|~ a^{20}=1,~b^{4}=1,~b^{-1}ab = a^ {r}\rangle$,~$r=3,13$;\\
  $\langle a,~b~| ~a^{20}=1,~b^{4}=a^{10},~b^{-1}ab = a^ {3}\rangle$;\\
  $\langle a,~b~|~ a^{21}=1,~b^{6}=1,~b^{-1}ab = a^ {r}\rangle$,~$r=2,10$;\\
  $\langle a,~b~|~ a^{28}=1,~b^{6}=a^{\ell},~b^{-1}ab = a^ {11}\rangle,~\ell=14,28$;\\
$\langle a,~b~|~ a^{30}=1,~b^{4}=1,~b^{-1}ab = a^ {17}\rangle$;\\
$\langle a,~b~| ~a^{36}=1,~b^{6}=a^{\ell},~b^{-1}ab = a^ {7}\rangle,~\ell=6,36$;\\
$\langle a,~b~|~ a^{42}=1,~b^{6}=1,~b^{-1}ab = a^ {r}\rangle$,~$r=11,19$;
\end{quote}

\noindent where $\lambda_{n}$ is a generator of $\mathcal{U}(\mathbb{Z}/n\mathbb{Z})$ and $\varphi$ denotes the Euler's phi function.
\end{theorem}
{\bf Proof.} Let $G$ be as in (\ref{Eq3}). Suppose $G$ has the {\sf cut}-property. Since the {\sf cut}-property is quotient closed, we have that the cyclic group $G/\langle a \rangle$ also has the {\sf cut}-property and hence,
\begin{equation}\label{e1}
  t = 2, 3, 4 ~{\rm or}~ 6.
\end{equation}
 Denote by $o_{n}(r)$, the order of $r$ modulo $n$. If $o_{n}(r) = d$, then one can check that $(\langle a,\, b^{d}\rangle,~ \langle  a^{\alpha}b^{d}\rangle)$ is a strong Shoda pair of $G$, for some $0\leq \alpha < n$. For instance, the following may be a choice of $\alpha$ for various cases:
 \begin{itemize}
   \item If $t=d$, then $\alpha=n-\ell$.
   \item If $t\neq d$ but $\ell=n$, then $\alpha=n$.
   \item If $t\neq d$, $\ell\neq n$ and $\operatorname{g.c.d.}(\frac{t}{d},n)$ divides $\ell$, where $\operatorname{g.c.d.}(\frac{t}{d},n)$ denotes the greatest common divisor of $\frac{t}{d}$ and $n$, then $\alpha$ is a solution of $\frac{t}{d}x\equiv-\ell~({\rm mod}~n)$.
   \item For all the other cases, $\alpha=\left\{
                 \begin{array}{ll}
                   \frac{m-\ell}{2}, & \hbox{if $\frac{t}{d}=2;$} \\
                  \frac{m-\ell}{3} ~{\rm or}~\frac{2m-\ell}{3}, & \hbox{if $\frac{t}{d}=3$,}
                 \end{array}
               \right.$  where $m$ denotes the greatest divisor of $n$, which is co-prime to $\frac{t}{d}$.

 \end{itemize}
   In view of (\ref{Eq2}), the centre of the simple component of $\mathbb{Q}[G]$ corresponding to this strong Shoda pair, computed using Proposition \ref{p11}, must be either $\mathbb{Q}$ or imaginary quadratic. Therefore, we must have \begin{equation}\nonumber \label{e2} \frac{\varphi(n)}{d}\leq 2.\end{equation} Consequently, $[\mathcal{U}(\mathbb{Z}/n\mathbb{Z}):\langle r \rangle]\leq 2$. We consider the following possible cases:

 \begin{description}
   \item[Case(I)] \underline{\textbf{$[\mathcal{U}(\mathbb{Z}/n\mathbb{Z}):\langle r \rangle]=1$}}\\

       \vspace{-0.25cm}
       Since $\mathcal{U}(\mathbb{Z}/n\mathbb{Z})=\langle r \rangle$, $o_{n}(r)$ divides $t$ and $t =2, 3, 4$ or $6$ (by (\ref{e1})), we have in this case $o_{n}(r)=\varphi(n)=2, 3, 4$ or $6$.
\begin{description}
  \item[Subcase (i)] \underline{$\varphi(n) = 2$}. Thus, $n= 3, 4, 6$ and $r=n-1$. In this case, for any~$n$, the admissible values of $t$ are $2, 4, 6$. Therefore, for suitable choices of $\ell$, $G =  \langle a,~b~|~ a^{n}=1,~b^{t}=a^{\ell},~b^{-1}ab = a^ {n-1}\rangle,$ where $n=3,4$ or $6$ and $t=2,4$ or $6.$ \\
      \vspace{-0.25cm}
  \item[Subcase (ii)] \underline{$\varphi(n) \neq 2$}. Observe that in this case, $t=o_{n}(r)=\varphi(n)=4$~or~$6$. As $\mathcal{U}(\mathbb{Z}/n\mathbb{Z})$ is cyclic, $n=5,7,9,10,14$ or $18$. Hence, for suitable choices of~$\ell$,
 $G \cong \langle a,~b~|~ a^{n}=1,~b^{t}=a^{\ell},~b^{-1}ab = a^ {\lambda_{n}}\rangle,$ where$~n=5,7,9,10,14$ or $18$ and $t=\varphi(n)$.
  \end{description}
 Note that if $n$ is even and $\ell$ is odd, then $t\neq4,\,6$. For otherwise, $(G,\langle a^{2}\rangle)$ is a strong Shoda pair of $G$ and centre of the simple component corresponding to this pair is $\mathbb{Q}(\zeta_{2t})=\mathbb{Q}(\zeta_{8})$ or $\mathbb{Q}(\zeta_{12})$. Therefore, we get that any group in this case, which has the {\sf cut}-property, must be one of the following:

 \begin{quote} $ \langle a,~b~|~ a^{n}=1,~b^{t}=1,~b^{-1}ab = a^ {n-1}\rangle,~~t=2,4,6,~~n=3,4,6;$ \\
 $\langle a,~b~| ~a^{4}=1,~b^{t}=a^{2},~b^{-1}ab = a^ {n-1}\rangle,~t=2,4,6$; \\
  $ \langle a,~b~|~ a^{6}=1,~b^{2}=a^{3},~b^{-1}ab = a^ {n-1}\rangle;$ \\
  $\langle a,~b~|~ a^{n}=1,~b^{\varphi(n)}=1,~b^{-1}ab = a^ {\lambda_{n}}\rangle,~n=5,7,9,10,14,18$.\end{quote}
  \item[Case(II)]\underline{\textbf{$[\mathcal{U}(\mathbb{Z}/n\mathbb{Z}):\langle r \rangle]=2$ and $\mathcal{U}(\mathbb{Z}/n\mathbb{Z})$ is cyclic}}  \\

      \vspace{-0.25cm} Clearly, $o_{n}(r)=\frac{\varphi(n)}{2}.$ As $o_{n}(r)$ divides $t$, $\varphi(n)=2,4,6,8$ or $12$. Since $G$ is non abelian, $\varphi(n)\neq 2$. Also, $\mathcal{U}(\mathbb{Z}/n\mathbb{Z})$ is cyclic yields following subcases:\begin{description}
                             \item[Subcase(i)] \underline{$\varphi(n)=4$, $n=5,10$.} In this case $t=2,~ 4$ or $6$; and $r=n-1$. Notice that the centre of the simple component corresponding to the strong Shoda pair $(\langle a, b^{2} \rangle, \langle b^{2} \rangle)$ is $\mathbb{Q}(\zeta_{n}+ \zeta_{n}^{-1}),$ which is a real quadratic extension of $\mathbb{Q}$. Hence, this case yields no group with the {\sf cut}-property.
      \item[Subcase(ii)] \underline{$\varphi(n)=6$, $n=7,9,14,18$.} In this case, $t=3$ or $6$, as $o_{n}(r)$ divides~$t$. As noted in Case(I), if $n=14$ and $t=6$, then $\ell \neq 7$ and hence, we must have $b^{t}=1$ i.e. $G \cong \langle a,~b~|~ a^{n}=1,~b^{\varphi(n)}=1,~b^{-1}ab = a^ {\lambda_{n}^{2}}\rangle$ or $G \cong \langle a,~b~|~ a^{n}=1,~b^{\frac{\varphi(n)}{2}}=a^{\ell},~b^{-1}ab = a^ {\lambda_{n}^{2}}\rangle, ~\ell=0~({\rm mod}~ 7)$, if $n=7,14$. Similar arguments yield that $\ell\equiv0~({\rm mod}~ 6)$, if $n=18$ and $t=6$. Otherwise, $\ell\equiv0~({\rm mod}~ 3)$, for $n=9,18$.
        Further, in view of \cite{Basma}, the groups listed above are isomorphic to one of the following:
  \begin{equation}\label{Eq5}
  \langle a,~b~|~ a^{n}=1,~b^{\frac{\varphi(n)}{j}}=1,~b^{-1}ab = a^ {\lambda_{n}^{2}}\rangle,~~j=1,2,~~n=7,9,14,18.   \end{equation}
                             \item[Subcase(iii)] \underline{$\varphi(n)=12$, $n=13,26$.} Here, $(\langle  a \rangle, \langle 1 \rangle)$ is a strong Shoda pair and the centre of the corresponding simple component is contained in $\mathbb{Q}(\zeta_{n}+ \zeta_{n}^{-1})$, which is a real extension of $\mathbb{Q}$. Hence, this case also yields no group with the {\sf cut}-property.
                           \end{description}

Therefore, if $[\mathcal{U}(\mathbb{Z}/n\mathbb{Z}):\langle r \rangle]= 2$ and $\mathcal{U}(\mathbb{Z}/n\mathbb{Z})$ is cyclic group generated by $\lambda_{n}$, then $G$ is isomorphic to one of the groups in (\ref{Eq5}).

 \item[Case(III)] \underline{\textbf{$[\mathcal{U}(\mathbb{Z}/n\mathbb{Z}):\langle r \rangle]=2$ and $\mathcal{U}(\mathbb{Z}/n\mathbb{Z})$ is non-cyclic}}\\
                              \vspace{-0.25cm}

                              For each group in this case, we enlist below the strong Shoda pair of $G$, \linebreak denoted by \emph{SSP}, which yields the simple component of $\mathbb{Q}[G]$, whose centre is different from $\mathbb{Q}$ or imaginary quadratic and hence restricts $G$ from having the {\sf cut}-property. For each $n$, we list the choice of parameters $t,\,r$ and $\ell$, which define the group $G$, as given in (\ref{Eq3}).
                              As $\mathcal{U}(\mathbb{Z}/n\mathbb{Z})$ is non-cyclic, $o_{n}(r)=\frac{\varphi(n)}{2}$ and  $o_{n}(r)$ divides $t$, we have the following subcases:

                              \begin{description}
                             \item[Subcase(i)] \underline{ $\varphi(n)=4$,~ $n=8,12,$~ $t=2,4$ or $6$.}\\

               \underline{$n=8$}~~~~
\begin{tabular}{|c|c|c|l|}
 \hline t & r & $\ell$ & \emph{SSP} \\ \hline
  2 & 7 & 4,8 & $(\langle a \rangle,~\langle 1 \rangle)$\\\hline
      4 & 3 & 4 & $(\langle a,\,b^{2} \rangle,~\langle a^{2}b^{2} \rangle)$\\\hline
   4& 5 & 2 & $(\langle a,\,b^{2} \rangle,~\langle a^{3}b^{2} \rangle)$\\\hline
   4& 7 & 8 & $(\langle a,\,b^{2} \rangle,~\langle b^{2} \rangle)$\\\hline
  6 & 3,7 & 8 & $(\langle a,\,b^{2} \rangle,~\langle b^{2} \rangle)$\\\hline
   6& 5 & 8 & $(\langle a,\,b \rangle,~\langle a^{4},\,b^{3} \rangle)$\\\hline

6&5  & 4 & $(\langle a,\,b^{2} \rangle,~\langle a^{4}b^{2} \rangle)$\\\hline                                              \end{tabular}\vspace{0.25cm}

                                              \underline{$n=12$}~~
\begin{tabular}{|c|c|c|l|}
 \hline t & r & $\ell$ & \emph{SSP} \\ \hline
  2 & 5 & 3 & $(\langle a,\,b \rangle,~\langle a^{4} \rangle)$\\\hline
   2& 11 & 6,12 & $(\langle a \rangle,~\langle 1 \rangle)$\\\hline

  4 &  5& 3 & $(\langle a,\,b \rangle,~\langle a^{2} \rangle)$\\\hline
 4&5  & 6 & $(\langle a,\,b \rangle,~\langle ab^{2} \rangle)$\\\hline
4& 7 & 2,12 & $(\langle a,\,b \rangle,~\langle a^{3} \rangle)$\\\hline
4& 11 & 12 & $(\langle a,\,b^{2} \rangle,~\langle b^{2} \rangle)$\\\hline
4& 11 & 6 & $(\langle a,\,b^{2} \rangle,~\langle a^{9}b^{2} \rangle)$\\\hline

6 &  5& 12 & $(\langle a,\,b \rangle,~\langle a^{4},\,b^{3} \rangle)$\\\hline
 6& 5 & 3 & $(\langle a,\,b \rangle,~\langle a^{2} \rangle)$\\\hline

  6& 7& 2,4 & $(\langle a,\,b \rangle,~\langle a^{3} \rangle)$\\\hline

  6& 11& 12 & $(\langle a,\,b^{2} \rangle,~\langle a^{4}b^{2} \rangle)$\\\hline
 6& 11 & 6 & $(\langle a,\,b^{2} \rangle,~\langle a^{10}b^{2} \rangle)$\\\hline
 \end{tabular}\vspace{0.25cm}

 \item[Subcase(ii)] \underline{$\varphi(n)=8$,~ $n=15,16,20,30,$~ $t=4$.}\\

                         \underline{$n=15$}~~
\begin{tabular}{|c|c|c|l|}
 \hline t & r & $\ell$ & \emph{SSP} \\ \hline
  4 & 7 & 15 & $(\langle a,\,b \rangle,~\langle a^{3} \rangle)$\\\hline
         \end{tabular}\vspace{0.25cm}

\underline{$n=20$}~~
\begin{tabular}{|c|c|c|l|}
 \hline t & r & $\ell$ & \emph{SSP} \\ \hline
  4 & 13 & 5 & $(\langle a,\,b \rangle,~\langle a^{2} \rangle)$\\\hline
  4       &13  & 10 & $(\langle a,\,b \rangle,~\langle ab^{2},\,a^{4} \rangle)$\\\hline \end{tabular}\vspace{0.25cm}

  \underline{$n=30$}~~
\begin{tabular}{|c|c|c|l|}
 \hline t & r & $\ell$ & \emph{SSP} \\ \hline
  4 & 7 & 5,30 & $(\langle a,\,b \rangle,~\langle a^{3} \rangle)$\\\hline
  4       &17 & 15 & $(\langle a,\,b \rangle,~\langle a^{2} \rangle)$\\\hline \end{tabular}\vspace{0.25cm}

                        \item[Subcase(iii)] \underline{ $\varphi(n)=12$,~ $n=21,28,36,42,$~ $t=6$.}\\

                        \underline{$n=21$}~~
\begin{tabular}{|c|c|c|l|}
 \hline t & r & $\ell$ & \emph{SSP} \\ \hline
  6 & 5 & 21 & $(\langle a \rangle,~\langle 1 \rangle)$\\\hline
    6&10  & 7 & $(\langle a,\,b \rangle,~\langle a^{15} \rangle)$\\\hline
                \end{tabular}\vspace{0.25cm}

                        \underline{$n=28$}~~
\begin{tabular}{|c|c|c|l|}
 \hline t & r & $\ell$ & \emph{SSP} \\ \hline
  6 & 3 & 14,28 & $(\langle a \rangle,~\langle 1 \rangle)$\\\hline
   6 & 5 & 28 & $(\langle a,\,b \rangle,~\langle a^{4},\,b^{3} \rangle)$\\\hline
    6&5  & 7 & $(\langle a,\,b \rangle,~\langle a^{2} \rangle)$\\\hline
                \end{tabular}\vspace{0.25cm}

\underline{$n=36$}~~
\begin{tabular}{|c|c|c|l|}
 \hline t & r & $\ell$ & \emph{SSP} \\ \hline
  6 & 5 & 36 & $(\langle a,\,b \rangle,~\langle a^{4},\,b^{3} \rangle)$\\\hline
    6&  5& 9 & $(\langle a,\,b \rangle,~\langle a^{2} \rangle)$\\\hline
    6& 11 & 36 & $(\langle a,\,b^{2} \rangle,~\langle a^{4}b^{2} \rangle)$\\\hline
    6&11  & 18 & $(\langle a,\,b^{2} \rangle,~\langle a^{10}b^{2} \rangle)$\\\hline
                \end{tabular}\vspace{0.25cm}
   \vspace{0.25cm}

       \underline{$n=42$}~~
\begin{tabular}{|c|c|c|l|}
 \hline t & r & $\ell$ & \emph{SSP} \\ \hline
  6 & 5 & 42 & $(\langle a \rangle,~\langle 1 \rangle)$\\\hline
   6 &5, 11 & 21 & $(\langle a,\,b \rangle,~\langle a^{2} \rangle)$\\\hline
6& 19 & 7,21 & $(\langle a,\,b \rangle,~\langle a^{2} \rangle)$\\\hline
    6& 19 & 14 & $(\langle a,\,b \rangle,~\langle a^{3} \rangle)$\\\hline
                   \end{tabular}\vspace{0.25cm}
                                                       \end{description}

              Clearly, the groups in subcases (i)-(iii) (with different set of parameters) which are isomorphic to one of the groups tabulated above, fail to have the {\sf cut}-property. The remaining ones are isomorphic to one of the following:
         \begin{quote}
  $\langle a,~b~|~ a^{8}=1,~b^{t}=1,~b^{-1}ab = a^ {r}\rangle,~t=2,4,~r=3,5$.\\
   $\langle a,~b~|~ a^{12}=1,~b^{t}=1,~b^{-1}ab = a^ {5}\rangle,~t=2,4$.\\
  $\langle a,~b~|~ a^{12}=1,~b^{t}=a^{\ell},~b^{-1}ab = a^ {7}\rangle,~t=2,6,~\ell=t,12$.\\
  $\langle a,~b~| ~a^{15}=1,~b^{4}=1,~b^{-1}ab = a^ {2}\rangle.$\\
  $\langle a,~b~|~ a^{16}=1,~b^{4}=1,~b^{-1}ab = a^ {r}\rangle$,~$r=3,5$.\\
  $\langle a,~b~| ~a^{20}=1,~b^{4}=1,~b^{-1}ab = a^ {r}\rangle$,~$r=3,13$.\\
  $\langle a,~b~| ~a^{20}=1,~b^{4}=a^{10},~b^{-1}ab = a^ {3}\rangle$.\\
  $\langle a,~b~|~ a^{21}=1,~b^{6}=1,~b^{-1}ab = a^ {r}\rangle$,~$r=2,10$.\\
  $\langle a,~b~| ~a^{28}=1,~b^{6}=a^{\ell},~b^{-1}ab = a^ {11}\rangle,~\ell=14,28$.\\
$\langle a,~b~|~ a^{30}=1,~b^{4}=1,~b^{-1}ab = a^ {17}\rangle$.\\
$\langle a,~b~|~ a^{36}=1,~b^{6}=a^{\ell},~b^{-1}ab = a^ {7}\rangle,~\ell=6,36$.\\
$\langle a,~b~|~ a^{42}=1,~b^{6}=1,~b^{-1}ab = a^ {r}\rangle$,~$r=11,19$.\end{quote}             \end{description}

For isomorphism of groups, we have again used the results of \cite{Basma}.\\
\vspace{-0.25cm}

It only remains to check that each group $G$ listed in the theorem has the {\sf cut}-property. Since metacyclic groups are normally monomial \cite{Bas}, and hence  strongly monomial \cite{BM1}, this can be done using Theorem 3.1 of \cite{jesper} (see also \cite{BM2}, Theorem 2) on the rank of the free abelian component of $\mathcal{Z}(\mathcal{U}(\mathbb{Z}[G]))$ and Theorem 1 of \cite{BM1} on the computation of strong Shoda pairs of $G$. $\Box$

\begin{remark}
In view of {\rm(\cite{hupp}, p.\,545)}, an alternative approach to determine metacyclic groups with the {\sf cut}-property would be to compute the character fields of their complex irreducible representations {\rm\cite{basm}}.
\end{remark}
\begin{remark}Recall that a group $G$ is said to have the normalizer property, if \linebreak $G.\mathcal{Z}(\mathcal{U}({\mathbb{Z}[G]}))=N_{\mathcal{U}}(G)$, the normalizer of $G$ in $\mathcal{U}({\mathbb{Z}[G]})$. It is known {\rm(\cite{Hai}, Lemma 2.1)} that the {\sf cut}-property implies the normalizer property. It is interesting to note that although it is easy to see that all metacyclic groups have the normalizer property, only a few of them have the cut property.
\end{remark}

 \subsection{The upper central series} It is well known (\cite{satya}, Theorem 2.6) that, for a finite group $G$, the central height of $V(\mathbb{Z}[G])$, the group of units of augmentation 1 in $\mathbb{Z}[G]$, is at most $2$. In this section, we give a classification of finite metacyclic groups $G$ for the permissible values of the central height of $V(\mathbb{Z}[G])$.
\begin{theorem}\label{Th6}   Let $G$ be a finite metacyclic group and let $V(\mathbb{Z}[G])$ be the group of units of augmentation $1$ in  $\mathbb{Z}[G]$.

\begin{description}
   \item[(I)] The  central height of  $V(\mathbb{Z}[G])$ equals  $0$
if, and only if, $G$ is isomorphic to one of following:
  \begin{quote}
               $\langle a,~b~|~ a^{3}=1,~b^{2}=1,~b^{-1}ab=a^{2}\rangle$\\
                 $\langle a,~b~|~ a^{5}=1,~b^{4}=1,~b^{-1}ab=a^{2}\rangle \\ $
                $\langle a,~b~|~ a^{7}=1,~b^{6}=1,~b^{-1}ab=a^{3}\rangle$\\
                 $\langle a,~b~|~ a^{9}=1,~b^{6}=1,~b^{-1}ab=a^{2}\rangle$\\
                $\langle a,~b~|~ a^{7}=1,~b^{3}=1,~b^{-1}ab=a^{2}\rangle$\\
                 $\langle a,~b~|~ a^{15}=1,~b^{4}=1,~b^{-1}ab=a^{2}\rangle$\\
                 $\langle a,~b~|~ a^{21}=1,~b^{6}=1,~b^{-1}ab=a^{2}\rangle.$

\end{quote}
\item[(II)] The central height of $V(\mathbb{Z}[G])$ equals $2$ if, and only if, $G$ is isomorphic to  $$\langle a,~b~|~ a^{n}=1,~b^{2}=a^{\frac{n}{2}},~b^{-1}ab=a^{-1}\rangle,$$ where $n$ is a positive integer divisible by $4$.
      \item[(III)]For all other groups $G$, the central height of $V(\mathbb{Z}[G])$ equals $1$.
      \end{description}
      \end{theorem}
\noindent{\bf  Proof.} Let $G$ be as in (\ref{Eq3}).
 \begin{description}
\item[(I)] Clearly, if the central height of $V(\mathbb{Z}[G])$ equals $0$, i.e., $V(\mathbb{Z}[G])=\langle 1 \rangle$, then $G$ must be one of the groups listed in Theorem \ref{Th7} with $\mathcal{Z}(G)=\langle 1 \rangle$. Observe that $1\neq a^{i}b^{j}\in \mathcal{Z}(G)$ if, and only if, for some $i,j$, not both $0$, $ri\equiv i({\rm mod}~n)$ and $r^{j}\equiv 1({\rm mod}~n)$, where $0\leq i < n$ and $0\leq j < t$. This occurs if, and only if, either $\operatorname{g.c.d.}(r-1,n)\neq 1$ or $b^{o_{n}(r)}\neq 1$. These observations yield (I).

  \item[(II)] In view of (\cite{saty}, Theorem 3.7), the central height of  $V(\mathbb{Z}[G])$ is $2$
   if, and only if, $G$ is a $Q^{*}$ group, i.e., there exists an abelian subgroup $H$ of $G$ of index $2$ and an element $x$ of order $4$ such that  $G = \langle H,\, x\rangle$ and the following holds: \begin{quote} (i) $x^{-1}hx = h^{-1}$ for each $h \in H$; \\(ii) $x^{2} = y^{2} $ for some $y \in H$; \\(iii) $H$ is not an elementary Abelian $2$-group.\end{quote} We thus examine which metacyclic groups are $Q^{*}$ groups. As the index of $H$ in $G$ is $2$, we have that $a^{2} \in H$, and therefore $H \cap \langle a \rangle = \langle a \rangle $ or $ \langle a^{2} \rangle.$  We take these two cases seperately.  \noindent
\vspace{.25cm}\\ \underline{$H \cap \langle a \rangle = \langle a \rangle $.} \\
The index of $H$ in $G$ being $2$, we must have that  $ H = \langle a,\, b^{2} \rangle.$  Furthermore, as $H$ is abelian, $b^{2}$ commutes with $a$ and hence $b^{2}\in \mathcal{Z}(G)$. This gives $b^{4}=1$ since $x^{-1}b^{2}x= b^{-2}$ (by (i)).
As $x \in Hb$, i.e., $x=h'b$, for some $h'\in H$, it follows from (i) that $b^{-1}hb=h^{-1}$, for every $h \in H$ and hence $ x^{2}=h'bh'b=b^{2}$. This implies that the order of $b$ must be $4$. Also by (ii), there exists an element $y \in H$ of order $4$ with $x^{2}= y^{2}$. However, $H$ has an element of order $4$ if, and only if, $4$ divides $n$, and in that case, the only elements in $H$ of order $4$ are $(a^{\frac{n}{4}})^{i}b^{2j}$, with $ i=1,\, 3$ and $j=0,\, 1$. Further, notice that the square of all of these four elements  is $a^{\frac{n}{2}}$. Consequently, we must have $b^{2}= a^{\frac{n}{2}}$ and $ G \cong \langle a, b~|~ a^{n}=1,~b^{2}=a^{\frac{n}{2}},~b^{-1}ab=a^{-1}\rangle $, where $n$ is divisible by 4.
\vspace{.25cm}\\
\underline{$H \cap \langle a \rangle = \langle a^{2} \rangle $.}  \\ In this case, we have  $ H = \langle a^{2},\, b \rangle$ or  $\langle a^{2},\, ab \rangle,$ depending upon whether \linebreak $b$ belongs to $H$ or not. Suppose  $ H = \langle a^{2},\, b \rangle.$ Since $H$ is abelian, we have that $a^{2}$ commutes with~$b$, so that $a^{2}\in \mathcal{Z}(G)$. It thus follows by (i) that $a^{4} =1$. However, $\langle a \rangle$ is a normal subgroup of non abelian group $G$, which yields that the order of $a$ cannot be $2$. Therefore, the order of $a$ is $4$ and $bab^{-1}$ must be $a^{-1}$. Further as $x\in Ha$, $x^{-1}bx=b^{-1}$ implies that $a^{-1}ba=b^{-1}$. Consequently, $a^{2}=b^{2}$ and $G \cong \langle a, b~|~ a^{4}=1,~b^{2}=a^{2},~b^{-1}ab=a^{-1}\rangle $. \par Similarly, if $H =  \langle a^{2}, \,ab \rangle,$ then $G \cong \langle a, b~|~ a^{4}=1,~b^{2}=a^{2},~b^{-1}ab=a^{-1}\rangle $. $\Box$
           \end{description}
\bibliographystyle{amsplain}
\bibliography{ReferencesBMP}

\providecommand{\bysame}{\leavevmode\hbox to3em{\hrulefill}\thinspace}
\providecommand{\MR}{\relax\ifhmode\unskip\space\fi MR }
\providecommand{\MRhref}[2]{%
  \href{http://www.ams.org/mathscinet-getitem?mr=#1}{#2}
}
\providecommand{\href}[2]{#2}
\begin{thebibliography}{10}

\bibitem{satya}
S.~R. Arora, A.~W. Hales, and I.~B.~S. Passi, \emph{Jordan decomposition and
  hypercentral units in integral group rings}, Comm. Algebra \textbf{21}
  (1993), no.~1, 25--35.

\bibitem{saty}
S.~R. Arora and I.~B.~S. Passi, \emph{Central height of the unit group of an
  integral group ring}, Comm. Algebra \textbf{21} (1993), no.~10, 3673--3683.

\bibitem{BM1}
G.K. Bakshi and S.~Maheshwary, \emph{The rational group algebra of a normally
  monomial group}, J. Pure Appl. Algebra \textbf{218} (2014), no.~9,
  1583--1593.

\bibitem{BM2}
\bysame, \emph{On the index of a free abelian subgroup in the group of central
  units of an integral group ring}, J. Algebra \textbf{434} (2015), 72--89.

\bibitem{Bas}
B.~G. Basmaji, \emph{Monomial representations and metabelian groups}, Nagoya
  Math. J. \textbf{35} (1969), 99--107.

\bibitem{Basma}
\bysame, \emph{On the ismorphisms of two metacyclic groups}, Proc. Amer. Math.
  Soc. \textbf{22} (1969), 175--182.

\bibitem{basm}
\bysame, \emph{Complex representations of metacyclic groups}, Amer. Math.
  Monthly \textbf{86} (1979), no.~1, 47--48.

\bibitem{dark}
R.~Dark and C.~M. Scoppola, \emph{On {C}amina groups of prime power order}, J.
  Algebra \textbf{181} (1996), no.~3, 787--802.

\bibitem{Feit}
W.~Feit and J.~G. Thompson, \emph{Solvability of groups of odd order}, Pacific
  J. Math. \textbf{13} (1963), 775--1029.

\bibitem{Ferraz}
R.~A. Ferraz and J.~J. Sim{\'o}n-P{\i}nero, \emph{Central units in metacyclic
  integral group rings}, Comm. Algebra \textbf{36} (2008), no.~10, 3708--3722.

\bibitem{good}
E.~G. Goodaire, E.~Jespers, and C.~Polcino~Milies, \emph{Alternative loop
  rings}, North-Holland Mathematics Studies, vol. 184, North-Holland Publishing
  Co., Amsterdam, 1996.

\bibitem{Hig}
G.~Higman, \emph{The units of group-rings}, Proc. London Math. Soc. (2)
  \textbf{46} (1940), 231--248.

\bibitem{hupp}
B.~Huppert, \emph{Endliche {G}ruppen. {I}}, Die Grundlehren der Mathematischen
  Wissenschaften, Band 134, Springer-Verlag, Berlin-New York, 1967.

\bibitem{del}
E.~Jespers and {\'A}.~del R{\'{\i}}o, \emph{Group ring groups}, Volume 1:
  Orders and Generic Constructions of Units, De Gruyter, Berlin-Boston, 2015.

\bibitem{jesper}
E.~Jespers, G.~Olteanu, {\'A}.~del R{\'{\i}}o, and I.~Van~Gelder, \emph{Group
  rings of finite strongly monomial groups: central units and primitive
  idempotents}, J. Algebra \textbf{387} (2013), 99--116.

\bibitem{jesp3}
\bysame, \emph{Central units of integral group rings}, Proc. Amer. Math Soc.
  \textbf{142} (2014), 2193--2209.

\bibitem{jesp2}
E.~Jespers and M.~M. Parmenter, \emph{Construction of central units in integral
  group rings of finite groups}, Proc. Amer. Math. Soc. \textbf{140} (2012),
  no.~1, 99--107.

\bibitem{parm}
E.~Jespers, M.~M. Parmenter, and S.~K. Sehgal, \emph{Central units of integral
  group rings of nilpotent groups}, Proc. Amer. Math. Soc. \textbf{124} (1996),
  no.~4, 1007--1012.

\bibitem{Li}
Y.~Li and M.~M. Parmenter, \emph{Central units of the integral group ring
  {${\bf Z}A_5$}}, Proc. Amer. Math. Soc. \textbf{125} (1997), no.~1, 61--65.

\bibitem{Hai}
Z.~Li and J.~Hai, \emph{The normalizer property for integral group rings of
  finite solvable {T}-groups}, J. Group Theory \textbf{15} (2012), no.~2,
  237--243.

\bibitem{MacD}
I.~D. Macdonald, \emph{Some {$p$}-groups of {F}robenius and extra-special
  type}, Israel J. Math. \textbf{40} (1981), no.~3-4, 350--364 (1982).

\bibitem{Oli}
A.~Olivieri, {\'A}.~del R{\'{\i}}o, and J.~J. Sim{\'o}n, \emph{On monomial
  characters and central idempotents of rational group algebras}, Comm. Algebra
  \textbf{32} (2004), no.~4, 1531--1550.

\bibitem{Pass}
D.~S. Passman, \emph{Infinite crossed products}, Pure and Applied Mathematics,
  vol. 135, Academic Press Inc., Boston, MA, 1989.

\bibitem{rei}
I.~Reiner, \emph{Maximal orders}, Academic Press [A subsidiary of Harcourt
  Brace Jovanovich, Publishers], London-New York, 1975, London Mathematical
  Society Monographs, No. 5.

\bibitem{sehg}
J.~Ritter and S.~K. Sehgal, \emph{Integral group rings with trivial central
  units}, Proc. Amer. Math. Soc. \textbf{108} (1990), no.~2, 327--329.

\end{thebibliography}
\end{document}